\RequirePackage{ifpdf}
\ifpdf 
\documentclass[pdftex]{sigma}
\else
\documentclass{sigma}
\fi

\newcommand{\R}{\mathbb{R}}
\newcommand{\s}{\mathbf{S}}
\newcommand{\J}{\mathbf{J}}
\newcommand{\N}{\mathbb{N}}
\newcommand{\D}{{\mathcal D}}
\newcommand{\Dbold}{{\mbox{\boldmath $\D$}}}
\newcommand{\nn}{{\mathbf{n},\mathbf{p}}}
\newcommand{\E}{{\mathcal E}}
\newcommand{\F}{{\mathcal F}}
\newcommand{\G}{{\mathcal G}}
\newcommand{\Y}{{\mathcal Y}}
\newcommand{\taubold}{{\mbox{\small\boldmath $\tau$}}}
\newcommand{\sigmabold}{{\mbox{\small\boldmath $\sigma$}}}
\newcommand{\alphabold}{{\mbox{\small\boldmath $\alpha$}}}
\newcommand{\q}{\mathbf{q}}
\renewcommand{\j}{\mathbf{j}}
\newcommand{\p}{\mathbf{p}}
\newcommand{\x}{\mathbf{x}}
\newcommand{\y}{\mathbf{y}}
\newcommand{\z}{\mathbf{z}}
\newcommand{\n}{\mathbf{n}}
\newcommand{\0}{\mathbf{0}}

\numberwithin{equation}{section}

\begin{document}

\allowdisplaybreaks

\renewcommand{\thefootnote}{$\star$}

\renewcommand{\PaperNumber}{091}

\FirstPageHeading

\ShortArticleName{External Ellipsoidal Harmonics for the Dunkl--Laplacian}

\ArticleName{External Ellipsoidal Harmonics\\ for the Dunkl--Laplacian\footnote{This paper is a contribution to the Special
Issue on Dunkl Operators and Related Topics. The full collection
is available at
\href{http://www.emis.de/journals/SIGMA/Dunkl_operators.html}{http://www.emis.de/journals/SIGMA/Dunkl\_{}operators.html}}}

\Author{Hans VOLKMER}

\AuthorNameForHeading{H. Volkmer}

\Address{Department of Mathematical Sciences,
University of Wisconsin-Milwaukee,\\
P. O. Box 413,
Milwaukee, WI 53201,
USA}

\Email{\href{mailto:volkmer@uwm.edu}{volkmer@uwm.edu}}

\URLaddress{\url{http://www.uwm.edu/~volkmer/}}

\ArticleDates{Received September 22, 2008, in f\/inal form December 18,
2008; Published online December 23, 2008}

\Abstract{The paper introduces external ellipsoidal and external
sphero-conal $h$-harmonics for the Dunkl--Laplacian. These external
$h$-harmonics admit integral representations,  and they are
connected by a formula of Niven's type. External $h$-harmonics in
the plane are expressed in terms of Jacobi polynomials
$P_n^{\alpha,\beta}$ and Jacobi's functions
$Q_n^{\alpha,\beta}$ of the second kind.}


\Classification{33C52; 35C10}

\section{Introduction}
In a previous article \cite{Vo} the author introduced (internal) ellipsoidal $h$-harmonics $F_\nn$ and sphero-conal
$h$-harmonics $G_\nn$ for the Dunkl--Laplacian. The functions $F_\nn(\x)$ and $G_\nn(\x)$ are polynomials
in the variables $\x=(x_0,x_1,\dots,x_k)$ and they are $h$-harmonic, that is, they
satisfy the dif\/ferential-dif\/ference equation
\begin{gather}\label{1:DL}
\Delta_h u:=\big(\D_0^2+\D_1^2+\dots+\D_k^2\big) u=0,
\end{gather}
where $\D_j$ is def\/ined by
\[ \D_j u(\x):=\frac{\partial u(\x)}{\partial x_j}+\alpha_j \frac{u(\x)-u(\sigma_j(\x))}{x_j} \]
and $\sigma_j$ is the ref\/lection at the $j$th coordinate plane.
The parameters $\alpha_j$ enter the weight function
\begin{gather}\label{1:h}
 h(\x):=|x_0|^{\alpha_0}|x_1|^{\alpha_1}\cdots|x_k|^{\alpha_k}.
\end{gather}
Following the book \cite{DX} by Dunkl and Xu and the paper \cite{Xu} by Xu
we will assume that $\alpha_j\ge 0$ for each $j=0,1,\dots,k$
although one would expect that the range of validity can be extended analytically to the domain
$\alpha_j>-\frac12$ for each $j$. We will also exclude the case $k=1$, $\alpha_0=\alpha_1=0$
because we want the constant $\mu$ def\/ined below in \eqref{1:mu} to be positive.

The parity vector $\p=(p_0,p_1,\dots,p_k)$ has components in $\{0,1\}$ and
indicates that $F_\nn$ and~$G_\nn$ have parity $\p$ which means that
they are sums of monomials $x_0^{p_0+2q_0}x_1^{p_1+2q_1}\cdots x_k^{p_k+2q_k}$ with $q_j\in\N_0=\{0,1,2,\dots\}$.
The vector $\n=(n_1,n_2,\dots,n_p)$ counts the zeros of the corres\-ponding
Stieltjes quasi-polynomials inside $k$ adjacent open intervals.
If we set $m=2|\n|+|\p|:=2(n_1+\cdots+n_k)+(p_0+\cdots+p_k)$ then
$G_\nn$ is a spherical $h$-harmonic of degree $m$, that is, it is a $h$-harmonic homogeneous
polynomial of degree $m$. Moreover, $F_\nn$
is a polynomial of total degree $m$ which, however, is usually not homogeneous.
The sphero-conal $h$-harmonics of f\/ixed degree $m$ form an orthogonal
basis in the linear space of spherical $h$-harmonics of degree $m$
with respect to the inner product
\begin{gather}\label{1:inner}
 \langle u,v\rangle =\int_{\s^k} h^2(\x) u(\x)v(\x)\,d\s(\x) ,
\end{gather}
where $\s^k$ denotes the unit sphere in $\R^{k+1}$.
Equation \eqref{1:inner} uses the standard surface integral over the sphere so that
$\int_{\s^k} 1d\s$ is the surface area of $\s^k$.

In the present paper the theory is extended to include external ellipsoidal $h$-harmonics $\F_\nn$ and
external sphero-conal $h$-harmonics $\G_\nn$. Throughout, we will use calligraphic letters to denote
the ``external'' functions. The functions $\G_\nn$ and $\F_\nn$ are still solutions of \eqref{1:DL}
but they are no longer polynomials. They have proper decay as $\|\x\|\to\infty$ which makes
them suitable
for solving problems concerning the exterior of ellipsoids or spheres, respectively.
The external sphero-conal $h$-harmonics are easy to def\/ine. They are given by
\begin{gather}\label{1:ext}
\G_\nn(\x):=\|\x\|^{-2\mu-2m} G_\nn(\x) ,\qquad m:=2|\n|+|\p|,
\end{gather}
where
\begin{gather}\label{1:mu}
\mu:=\alpha_0+\alpha_1+\cdots+\alpha_k+\frac{k-1}{2} .
\end{gather}
The def\/inition of external ellipsoidal $h$-harmonics is postponed to Section~\ref{E}.

After introducing the fundamental solution of \eqref{1:DL} in Section~\ref{F}
we obtain integral representations of $\F_\nn$ and $\G_\nn$ in terms of
the fundamental solution in Section~\ref{I}.
Section~\ref{H} states some integral formulas related to Section~5.2 in~\cite{DX}.
As a main result of this paper, in Section~\ref{N}
we prove a Niven type formula that expresses $\F_\nn$ in terms of $\G_\nn$.
In the classical situation $k=2$, $\alphabold=0$, this formula was f\/irst given by Niven~\cite{N}
in an amazing award winning paper. Modif\/ied versions of some of
Niven's result are contained in the books by
Whittaker and Watson \cite{WW} and by Hobson \cite{Hob}.
It was the attempt to extend Niven's results to external $h$-harmonics
that led to this paper.
We also refer to paper \cite{KM} by Kalnins and Miller
which contains a generalization of Niven's work from a dif\/ferent perspective.

In Section \ref{P} we illustrate the results in the relatively simple but nontrivial planar case $k=1$.
In this situation the internal ellipsoidal and sphero-conal $h$-harmonics are expressed in terms of
Jacobi polynomials $P_n^{(\alpha,\beta)}$ while the external ellipsoidal $h$-harmonic $\F_\nn$ involves the
second solution $Q_n^{(\alpha,\beta)}$ of the dif\/ferential equation whose f\/irst solution is $P_n^{(\alpha,\beta)}$.

We should mention that there is an increasing number of papers in Applied Mathematics
that use ellipsoidal harmonics; for example, see \cite{B,D,DK}.
For applications of Dunkl operators and the associated Laplacian to mathematical physics,
special functions, probability theory and geometry, we refer to \cite{BO, DJ, Du2, Hi, R, RV}.

\section{External ellipsoidal harmonics}\label{E}
We consider the Fuchsian dif\/ferential equation
\begin{gather}\label{E:Fuchs}
\prod_{j=0}^k (t-a_j)\left[ v''+\sum_{j=0}^k \frac{\alpha_j+\frac12}
{t-a_j} v'\right] +\left[-\frac12\sum_{j=0}^k \frac{p_j\alpha_j A_j}{t-a_j}
+\sum_{i=0}^{k-1}\lambda_i t^i\right]v =0
\end{gather}
for the function $v(t)$ where the prime denotes dif\/ferentiation with respect to~$t$,
and
\begin{gather*}
 A_j:=\prod_{i=0 \atop i\ne j}^k (a_j-a_i).
\end{gather*}
The parameters $\alpha_j$ are from \eqref{1:h}, $\p=(p_0,p_1\dots,p_k)$ is a f\/ixed
parity vector,
\begin{gather}\label{E:a}
a_0<a_1<\dots<a_k
\end{gather}
are chosen numbers and $\lambda_0,\lambda_1,\dots,\lambda_{k-1}$ are real spectral parameters.

Let $\n=(n_1,n_2,\dots,n_k)\in\N_0^k$.
In \cite[Section~2]{Vo} we considered the quasi-polynomial $E_\nn$
originally introduced by Stieltjes \cite{St} in the special case $\p=0$. The quasi-polynomial $E_\nn$
has~$n_j$ zeros in $(a_{j-1},a_j)$ for each $j=1,2,\dots,k$, and it is
a solution of \eqref{E:Fuchs} for special values of the parameters
$\lambda_0,\lambda_1,\dots,\lambda_{k-1}$.
For this set of parameters, $-\frac{m}{2}$ and $\mu+\frac{m}{2}$ are the
exponents of equation \eqref{E:Fuchs} at inf\/inity, where $m$ and $\mu$ are according to
\eqref{1:ext}, \eqref{1:mu}.
Now $E_\nn$ is the Frobenius solution belonging to the exponent $-\frac{m}{2}$.
We introduce a second Frobenius solution~$\E_\nn$ of equation~\eqref{E:Fuchs}
belonging to the exponent $\mu+\frac{m}{2}$ at inf\/inity.
This solution is def\/ined for $t>a_k$ and it is normalized by the condition
that
\begin{gather*}
\lim_{t\to\infty} t^{\mu+\frac{m}{2}}\E_\nn(t)=1.
\end{gather*}
This normalization and the corresponding one of $E_\nn$ leads to the Wronskian
\begin{gather}\label{E:W}
\E_\nn(t)E'_\nn(t)-E_\nn(t)\E'_\nn(t)=(\mu+m)\prod_{j=0}^k (t-a_j)^{-\alpha_j-\frac12}
\end{gather}
which holds for $t>a_k$.
In the classical case $k=2$, $\alphabold=\0$, the second solution
$\E_\nn$ appears in \cite[\S~9.7]{A} and \cite[\S~100]{He1} .

We introduce ellipsoidal coordinates involving the
parameters \eqref{E:a}.
For every $(x_0,\dots,x_k)$ in the positive cone of $\R^{k+1}$
\begin{gather}\label{E:cone}
x_0>0,\ \dots,\ x_k>0,
\end{gather}
its ellipsoidal coordinates $t_0,t_1,\dots,t_k$ lie in the intervals
\begin{gather}\label{E:cube}
 a_k<t_0<\infty,\qquad a_{i-1}<t_i<a_i,\qquad i=1,\dots,k,
\end{gather}
and satisfy
\begin{gather*}
\sum_{j=0}^k \frac{x_j^2}{t_i-a_j} =1\qquad\text{for} \quad i=0,1,\dots,k .
\end{gather*}
These coordinates provide a bijective mapping between the positive
cone~\eqref{E:cone} and the cube~\eqref{E:cube}.

As explained in \cite[Section~3]{Vo}, equation
\eqref{1:DL} can be solved by the method of separation of variables in ellipsoidal coordinates.
In \cite{Vo} we considered the internal ellipsoidal $h$-harmonic
\[ F_\nn(\x)=\prod_{j=0}^k E_\nn(t_j), \]
where $(t_0,t_1,\dots,t_k)$ are ellipsoidal coordinates of $\x$.
We now introduce the external ellipsoidal $h$-harmonic $\F_\nn$ by
\begin{gather*}
\F_\nn(\x):=\E_\nn(t_0)\prod_{j=1}^k E_\nn(t_j)=\frac{\E_\nn(t_0)}{E_\nn(t_0)} F_\nn(\x).
\end{gather*}
This function is originally def\/ined in the positive cone \eqref{E:cone}
and can then be extended analytically to the complement of the
degenerate ellipsoid
\begin{gather}\label{E:deg}
 \sum_{j=0}^{k-1} \frac{x_j^2}{a_k-a_j}\le 1,\qquad x_k=0.
\end{gather}
We see this as follows.
The ellipsoidal coordinate $t_0>a_k$ is an analytic function of $\x$ for $\x$
outside the degenerate ellipsoid \eqref{E:deg}.
Further, $E_\nn$ and $\E_\nn$ are analytic functions
on the interval $(a_k,\infty)$, $E_\nn$ has no zero there and
$F_\nn(\x)$ is a polynomial in $\x$.

We note that $\F_\nn$ has parity $\p$, that is,
\[ \F_\nn(\sigma_j(\x))=(-1)^{p_j}\F_\nn(\x).\]
Therefore, by construction, $\F_\nn$ satisf\/ies equation \eqref{1:DL}.
It may be appropriate here to remark that a solution of \eqref{1:DL} is always
implicitly assumed to be def\/ined on an open set which is
invariant under the ref\/lections $\sigma_j$, $j=0,1,\dots,k$.
Of course, our domain of def\/inition of $\F_\nn$ has this property.

\section{The fundamental solution}\label{F}

We introduce the function
\begin{gather}\label{F:fund}
\Phi(\x,\z):=\frac{\Gamma(\mu)}{4\prod_{j=0}^k\Gamma(\alpha_j+\frac12)}
\int\limits_{[-1,1]^{k+1}}\!\!\!\!(\Psi(\x,\z,\taubold))^{-\mu}
\prod_{j=0}^k c_{\alpha_j}(1+\tau_j)(1-\tau_j^2)^{\alpha_j-1}\,d\taubold ,
\end{gather}
where
\begin{gather*}
 \Psi(\x,\z,\taubold):=\|\x\|^2-2(\tau_0x_0z_0+\dots+\tau_kx_kz_k)+\|\z\|^2\ge 0,
\end{gather*}
$c_\nu:=(B(\frac12,\nu))^{-1}$ with $B$ denoting the Beta function,
and $\mu$ is from \eqref{1:mu}.
If $\alpha_i=0$ then delete the factor containing $c_{\alpha_i}$,
set $\tau_i=1$ in $\Psi$ and omit the integration on $\tau_i$. If
$\alpha_j=0$ for all $j$ then $\Phi$ reduces to
the standard fundamental solution of Laplace's equation \cite[page 22]{E}.

If $|x_j|\ne |z_j|$ for at least one $j$, then $\Psi(\x,\z,\cdot)$ has a positive lower
bound on $[-1,1]^{k+1}$.
Therefore, $\Phi(\x,\z)$ is analytic everywhere in $\R^{k+1}\times\R^{k+1}$
except at points where $|x_j|=|z_j|$ for all $j$.
We recognize that $\Phi(\x,\z)=\Phi(\z,\x)$,
$\Phi(\gamma \x,\gamma \z)=|\gamma|^{-2\mu}\Phi(\x,\z)$ and $\Phi(\x,\z)>0$.
Moreover,
\begin{gather*}
\Phi(\0,\z)=\frac{\Gamma(\mu)}{4\prod\limits_{j=0}^k\Gamma(\alpha_j+\frac12)}\,\|\z\|^{-2\mu}.
\end{gather*}
If $\z\ne\0$ we can expand $\Phi(\cdot,\z)$ in a Taylor series at $\x=0$. Its coef\/f\/icients are given by the
following lemma.
\begin{lemma}\label{F:l1}
If $\z\ne\0$ then, for every $(q_0,q_1,\dots,q_k)\in\N_0^{k+1}$,
\begin{gather}
\left.\D_{x_0}^{q_0}\cdots\D_{x_k}^{q_k}\big(x_0^{q_0}\cdots x_k^{q_k}\big)\frac{\partial_{x_0}^{q_0}\cdots\partial_{x_k}^{q_k}\Phi(\x,\z)}
{q_0!\cdots q_k!}\right|_{\x=\0}\nonumber\\
\qquad{}=\D_{x_0}^{q_0}\cdots\D_{x_k}^{q_k}\Phi(\x,\z)\left.\right|_{\x=\0}
=(-1)^{|\q|}\D^{q_0}_{z_0}\cdots\D_{z_k}^{q_k}\Phi(\0,\z) ,\label{F:taylor}
\end{gather}
where $|\q|=q_0+\cdots+q_k$.
\end{lemma}
\begin{proof}
We note that the f\/irst equality in \eqref{F:taylor} is trivial so it is enough to prove the second one.
The proof of the second equality involves some lengthy formulas, so we will use
some common abbreviations. If
$\j=(j_0,j_1,\dots,j_k)$ is a multi-index we set $\j!=j_0!\cdots j_k!$,
$\partial_\x^\j=\partial_{x_0}^{j_0}\cdots \partial_{x_k}^{j_k}$
and $\x^\j=x_0^{j_0}\cdots x_k^{j_k}$.
By Faa di Bruno's formula (or by equation \cite[(5.4)]{Vo}), we write
\[
\partial_\x^\q f(\Psi(\x,\z,\taubold))
=\sum_{\j} 2^{|\q|-2|\j|}f^{(|\q|-|\j|)}
(\Psi(\x,\z,\taubold))
\frac{1}{\j!}\,\partial_\x^{2\j}(x_0-\tau_0z_0)^{q_0}
\cdots(x_k-\tau_kz_k)^{q_k},
\]
where $f(u):=u^{-\mu}$, and the summation is over all multi-indices $\j=(j_0,j_1,\dots,j_k)$
with $j_i\ge 0$ and $q_i-2j_i\ge 0$, $i=0,1,\dots,k$.
Computing the partial derivatives on the right-hand side and setting $\x=\0$ gives
\begin{gather*}
\left.\partial_\x^\q f(\Psi(\x,\z,\taubold))\right|_{\x=\0}\\
\qquad{}
=(-1)^{|\q|}\sum_{\j} 2^{|\q|-2|\j|}f^{(|\q|-|\j|)}
(\|\z\|^2)
\frac{\q!}{\j!(\q-2\j)!} (\tau_0z_0)^{q_0-2j_0}\cdots(\tau_kz_k)^{q_k-2j_k}
.
\end{gather*}
We now carry out the  integrations indicated in \eqref{F:fund} taking into account that
\[ c_{\alpha_j}\int_{-1}^1 \tau_j^m(1+\tau_j)(1-\tau_j^2)^{\alpha_j-1}\,d\tau_j=
\frac{m!}{\D_{z_j}^{m}(z_j^m)}.\]
Note that $\D_{z_j}^{m}(z_j^m)$ is a constant depending on $j$ and $m$ (the
``Dunkl factorial''.)
We obtain
\[
\left.\partial_\x^\q\tilde\Phi(\x,\z)\right|_{\x=\0}
=(-1)^{|\q|}\sum_{\j} 2^{|\q|-2|\j|}f^{(|\q|-|\j|)}(\|\z\|^2)
\frac{\q!}{\j!\,\D_\z^{\q-2\j}\,\z^{\q-2\j}}\,\z^{\q-2\j},
\]
where $\tilde\Phi$ is def\/ined by the right-hand side of \eqref{F:fund} but with
the normalization factor in front of the integral omitted.
Now we use the f\/irst identity in \eqref{F:taylor} to obtain
\begin{gather*}
 \left.\D_\x^\q\tilde\Phi(\x,\z)\right|_{\x=\0}
=(-1)^{|\q|}
\sum_{\j} 2^{|\q|-2|\j|}f^{(|\q|-|\j|)}(\|\z\|^2)
\frac{1}{\j!}\,\frac{\D_\z^\q\,\z^\q}{\D_\z^{\q-2\j}\,\z^{\q-2\j}}\z^{\q-2\j}\\
\phantom{\left.\D_\x^\q\tilde\Phi(\x,\z)\right|_{\x=\0}}{} =(-1)^{|\q|}
\sum_{\j} 2^{|\q|-2|\j|}f^{(|\q|-|\j|)}(\|\z\|^2)
\frac{1}{\j!}\,\D_\z^{2\j}\, \z^\q.
\end{gather*}
Finally, we apply \cite[Lemma~2]{Vo}  to the function  $A(v_0,\dots,v_k)=f(v_0+\cdots+v_k)$ and
obtain the second identity in~\eqref{F:taylor}.
\end{proof}

Since $\|\z\|^{-2\mu}$ is $h$-harmonic,
Lemma \ref{F:l1} shows that $\Phi(\x,\cdot)$
is $h$-harmonic for f\/ixed $\x$ and so, by symmetry,
$\Phi(\cdot,\z)$ is $h$-harmonic for f\/ixed $\z$.

Knowledge of the exact nature of the singularities of $\Phi$ will be immaterial
in this paper.
An analysis of the singularities would be along the following lines.
We see that the singularity of $\Phi(\x,\z)$ at $\x=\z$ is determined by integration over values
$\tau_j$ close to $1$ so it is convenient to substitute $\sigma_j=1-\tau_j$.
If all $z_j$ are nonzero we have the
elementary inequality
\[ \Psi(\x,\z,\taubold) \ge c(\|\x-\z\|^2+\|\sigmabold\|)\]
if $\x$ is close enough to $\z$ and $c$ is a positive constant; see the proof of
\cite[Theorem 5.5.7]{DX} for similar estimates.
It follows that $\Phi(\cdot,\z)$ and its normal derivative
on spheres $\|\x-\z\|=\epsilon$ do not grow faster as $\epsilon\to 0$ than in the classical
case $\alphabold=\0$.
Obvious modif\/ications of the growth rate are obtained when one or several $z_j$ vanish.
If $|x_j|=|z_j|$ for all $j$ but $x_i\ne z_i$ for at least one $i$ then the factor
$1+\tau_i$ appearing in \eqref{F:fund} comes into play and makes the singularity milder than
that at $\x=\z$.

We now express $\Phi$ in terms of the reproducing kernel $P_m(\x,\y)$
of the linear space of spherical $h$-harmonics of degree $m$ with respect to the inner product \eqref{1:inner}; see
\cite[Section 5.3]{DX}.
We normalize the reproducing kernel so that
\[ Y(\x)=\int_{\s^k}h^2(\x) P_m(\x,\y)Y(\y)\,d\s(\y) \]
for every spherical $h$-harmonic $Y$ of degree $m$ and every $\x\in\s^k$.
The generating function for Gegenbauer polynomials \cite[Def\/inition 1.4.10]{DX} is given by
\begin{gather}\label{F:gen}
(1-2us+s^2)^{-\mu}=\sum_{m=0}^\infty C_m^{\mu}(u)s^m .
\end{gather}
The convergence is uniform for $-1\le u\le 1$, $|s|\le q<1$.
In \eqref{F:gen} we set $s=\|\x\|/\|\z\|$ and $u=\|\x\|^{-1}\|\z\|^{-1}(\tau_0x_0z_0+\dots+\tau_kx_kz_k)$.
Then \eqref{F:fund} and the integral representation of the reproducing kernel
\cite[Theorem 5.5.5]{DX} yields
the following theorem.
\begin{theorem}\label{F:t1}
If $\|\x\|<\|\z\|$ then
\begin{gather}\label{F:Laplace}
\Phi(\x,\z)=\sum_{m=0}^\infty \frac{1}{2(\mu+m)} \|\x\|^m\|\z\|^{-2\mu-m}
P_m\left(\frac{\x}{\|\x\|},\frac{\z}{\|\z\|}\right) .
\end{gather}
The convergence is uniform for $\frac{\|\x\|}{\|\z\|}\le q<1$.
\end{theorem}

When $k=2$, $\alphabold=\0$, formula \eqref{F:Laplace} shows how to
expand the reciprocal distance between $\z$ and $\x$ into a series of Legendre polynomials.
This expansion is due to Laplace
and it marks the very beginning of the history of spherical harmonics \cite[page 3]{He1}.

\section[Integral formulas for external $h$-harmonics]{Integral formulas for external $\boldsymbol{h}$-harmonics}\label{I}

Multiplying \eqref{F:Laplace} by a spherical $h$-harmonic and integrating, we
obtain an integral representation of external spherical $h$-harmonics.

\begin{theorem}\label{I:t1}
Let $Y$ be a spherical $h$-harmonic of degree $m$, and let $\Y$ denote the corresponding
external spherical $h$-harmonic $\Y(\x):=\|\x\|^{-2\mu-2m}Y(\x)$. Then, for $\|\z\|>1$,
\begin{gather}\label{I:int1}
\Y(\z)=2(\mu+m)\int_{\s^k} h^2(\x)\Phi(\x,\z)Y(\x)\,
d\s(\x).
\end{gather}
\end{theorem}
In particular, Theorem \ref{I:t1} can be applied to a sphero-conal $h$-harmonic $Y=G_\nn$.

Let $Y$ be any spherical $h$-harmonic of degree $m$.
If $\|\z\|>1$ we use Green's formula for $h$-harmonic functions \cite[Lemma 5.1.5]{DX}
inside the unit ball
to obtain
\begin{gather}\label{I:eq1}
 0=\int_{\s^k} h^2(\x)\left(Y(\x)\frac{\partial \Phi(\x,\z)}{
\partial \nu_\x}-\Phi(\x,\z)\frac{\partial Y(\x)}{\partial \nu}\right)\,d\s(\x),
\end{gather}
where $\nu$ denotes the outward normal.
On the unit sphere we have
\[ \frac{\partial Y}{\partial \nu}=mY=m\Y,\qquad\frac{\partial \Y}{\partial \nu}=
(-2\mu-m)Y. \]
Therefore, adding \eqref{I:eq1} to \eqref{I:int1} gives
\begin{gather*}
\Y(\z)=\int_{\s^k} h^2(\x)\left(\Y(\x)\frac{\partial \Phi(\x,\z)}{\partial \nu_\x}-
\Phi(\x,\z)\frac{\partial\Y(\x)}{\partial \nu}\right)d\s(\x) .
\end{gather*}
When $\alphabold=0$ this formula corresponds exactly to the well-known representation
formula of harmonic
functions in terms of a fundamental solution \cite[page 34]{E}.
By using Green's formula in the region between two spheres
we see that we also have
\begin{gather}\label{I:rep2}
\Y(\z)=\int_{\s^k_R} h^2(\x)\left(\Y(\x)\frac{\partial \Phi(\x,\z)}{\partial \nu_\x}-
\Phi(\x,\z)\frac{\partial\Y(\x)}{\partial \nu}\right)d\s(\x)
\end{gather}
when $\|\z\|>R$, where
\[ \s^k_R:=\{\x: \|\x\|=R\}.\]

We introduce ellipsoidal coordinates by choosing numbers
$a_j$ according to \eqref{E:a},
and prove an integral representation of external ellipsoidal $h$-harmonics.

\begin{theorem}\label{I:t2}
Let $\J$ denote the ellipsoid
\begin{gather}\label{I:J}
\J=\left\{\y: \sum_{j=0}^k \frac{y_j^2}{t-a_j}=1\right\}
\end{gather}
for some fixed $t>a_k$.
Then, for all $\z$ exterior to $\J$,
\begin{gather}\label{I:int2}
 \F_\nn(\z)=\frac{2(\mu+m)}{E^2_\nn(t)}\prod_{j=0}^k (t-a_j)^{-\alpha_j-\frac12}
\int_\J h^2(\y)w(\y)\Phi(\y,\z)F_\nn(\y)\,d\J(\y) ,
\end{gather}
where $m=2|\n|+|\p|$ and $w$ is defined by
\begin{gather}\label{I:w}
 \frac{1}{w(\y)}:=\sum_{j=0}^k \frac{y_j^2}{(t-a_j)^2} .
\end{gather}
\end{theorem}
\begin{proof}
We will suppress the subscripts $\nn$ in this proof.
Both sides of equation \eqref{I:int2} are analytic in the exterior of $\J$ so
it will be suf\/f\/icient to prove the equation for $\|\z\|>R$,
where $R$ is chosen so large that
$\J$ f\/its inside the sphere $\s^k_R$.

We claim that formula  \eqref{I:rep2} is also true if we replace $\Y$ by $\F$.
To see this we expand $\F$ on the sphere $\s^k_R$ in a series
\[ \F(\x)=\sum_{q=0}^\infty Y_q(\x), \]
where $Y_q$ is a spherical $h$-harmonic of degree $q$.
Using an uniqueness theorem for $h$-harmonic functions def\/ined in the exterior of spheres
with proper decay at inf\/inity (the analogue of \cite[Theorem 2.9]{Du1})
we obtain
\[ \F(\x)=\sum_{q=0}^\infty R^{2\mu+2q} \|\x\|^{-2\mu-2q}Y_q(\x). \]
This series and the series obtained by term-by-term dif\/ferentiation
converge uniformly
to the proper limits on and exterior to $\s^k_R$. Therefore,
by taking limits we obtain
\[ \F(\z)=\int_{\s^k_R} h^2(\x)\left(\F(\x)\frac{\partial \Phi(\x,\z)}{\partial \nu_\x}-
\Phi(\x,\z)\frac{\partial\F(\x)}{\partial \nu}\right)d\s(\x) .\]
Using Green's formula in the region between $\J$ and $\s_R^k$, we also have
\begin{gather}\label{I:rep}
 \F(\z)=\int_\J h^2(\y)\left(
\F(\y)\frac{\partial\Phi(\y,\z)}{\partial \nu_\y}
-\Phi(\y,\z)\frac{\partial \F(\y)}{\partial\nu}
\right)\,d\J(\y),
\end{gather}
where $\nu$ denotes the outward unit normal vector.
We now apply Green's formula to the internal ellipsoidal $h$-harmonic
$F$ in the interior of $\J$ and obtain
\begin{gather}\label{I:Green}
\int_\J h^2(\y)\left(F(\y)\frac{\partial\Phi(\y,\z)}{\partial \nu_\y}
-\Phi(\y,\z)\frac{\partial F(\y)}{\partial\nu}
\right)\,d\J(\y)=0.
\end{gather}
We multiply \eqref{I:rep} by $E(t)$, \eqref{I:Green} by $\E(t)$
and subtract noting that
$E(t)\F(\y)=\E(t) F(\y)$ on $\J$.
We obtain
\[ E(t)\F(\z)=\int_\J h^2(\y)\Phi(\y,\z)\left(
\E(t)\frac{\partial F(\y)}{\partial\nu}
-E(t)\frac{\partial \F(\y)}{\partial\nu}
\right)\,d\J(\y). \]
Since the ellipsoidal coordinates are orthogonal, the
normal derivative is expressible in terms of the derivative
with respect to the ellipsoidal coordinate $t_0$:
\[ \frac{\partial}{\partial\nu}=2w(\y)\frac{\partial}{\partial t_0}.\]
Therefore, we f\/ind
\[ \F(\z)=2\frac{\E(t)E'(t)-E(t)\E'(t)}{E^2(t) }\int_\J h^2(\y)w(\y)\Phi(\y,\z)F(\y)\,d\J(\y).\]
We replace the Wronskian by \eqref{E:W} and \eqref{I:int2} is established.
\end{proof}
In Section \ref{N} the integral in \eqref{I:int2} will be evaluated
using the formulas given in Section \ref{H}.

Theorem \ref{I:t2} leads to an expansion of $\Phi$ in ellipsoidal $h$-harmonics.

\begin{theorem}\label{I:t3}
Let the ellipsoidal $t_0$-coordinate of $\z$ be larger than the corresponding one for $\y$.
Then
\begin{gather}\label{I:Heine}
 \Phi(\y,\z)=\sum_{\n,\p}\frac{e_\nn^2}{2(\mu+2|\n|+|\p|)} F_\nn(\y)\F_\nn(\z) ,
\end{gather}
where the positive constants $e_\nn$ are determined by
\[ e_\nn^2 \int_{\s^k} h^2(\x)G^2_\nn(\x)\,d\s(\x)=1 .\]
\end{theorem}
\begin{proof}
Let $\J$ be the ellipsoid \eqref{I:J} that contains $\y$.
The function $\Phi(\cdot,\z)$ is analytic and $h$-harmonic on and inside $\J$.
Therefore, it can be expanded  in ellipsoidal $h$-harmonics as in \cite[Section 7]{Vo}
and the expansion coef\/f\/icients can be evaluated using \eqref{I:int2}.
This gives \eqref{I:Heine} after a~simple calculation.
\end{proof}

When $k=2$, $\alphabold=\0$, formula \eqref{I:Heine}
can be found in \cite[page 172]{He2}.

\section[Integral formulas for spherical $h$-harmonics]{Integral formulas for spherical $\boldsymbol{h}$-harmonics}\label{H}

\begin{theorem}\label{H:t1}
Let $f$ be a homogeneous polynomial of degree $\ell$, and let $Y$ be a spherical
$h$-harmonic of degree $m$. If $\ell=m+2r$, $r=0,1,2,\dots$, then
\begin{gather}\label{H:eq1}
\int_{\s^k} h^2(\x)f(\x)Y(\x)\,d\s(\x)=
\frac{1}{2^{\ell-1}r!}\frac{\prod\limits_{i=0}^k \Gamma(\alpha_i+\frac12)}{\Gamma(m+\mu+r+1)}
\,\Delta_h^rY(\Dbold)f(\x),
\end{gather}
where $Y(\Dbold)=Y(\D_0,\dots,\D_k)$.
If $\ell<m$ or $\ell-m$ is odd then the integral in \eqref{H:eq1} vanishes.
\end{theorem}
\begin{proof}
Theorem 5.1.15 in \cite{DX} gives the expansion of $f$ in
spherical $h$-harmonics on $\s_k$. It follows from this expansion and
\cite[Theorem 5.1.6]{DX}
that the integral in \eqref{H:eq1}
vanishes when $\ell<m$ or when $\ell-m$ is odd.
If $\ell=m+2r$, $r=0,1,2,\dots$, then \cite[Theorem 5.1.15]{DX} shows that
\begin{gather}\label{H:eq2}
\int_{\s^k} h^2(\x)f(\x)Y(\x)\,d\s(\x)=\int_{\s^k} h^2(\x)\tilde{f}(\x)Y(\x)\,d\s(\x),
\end{gather}
where $\tilde f$ is a homogeneous polynomial of degree $m$ given by
\[ \tilde{f}(\x)=\frac{1}{4^rr!(m+\mu+1)_r}
\Delta_h^r f(\x). \]
We evaluate the integral on the right-hand side of \eqref{H:eq2} by
applying Theorem 5.2.4 in \cite{DX}.
We obtain
\[ \int_{\s^k} h^2(\x)f(\x)Y(\x)\,d\s(\x)=\frac{1}{c'2^m4^rr!(\mu+1)_m(m+\mu+1)_r}
\Delta_h^r Y(\Dbold) f(\x),\]
where
\begin{gather*}
 \frac{1}{c'}=\int_{\s^k} h^2(\x)\,d\s(x)=
\frac{2\prod\limits_{i=0}^k \Gamma(\alpha_i+\frac12)}{\Gamma(\mu+1)}.
\end{gather*}
After simplif\/ication of the constant factors we arrive at the desired equation \eqref{H:eq1}.
\end{proof}

In the special case $k=2$, $\alphabold=\0$ formula \eqref{H:eq1} is due to Hobson
\cite[\S~100]{Hob}.

We state a consequence of Theorem \ref{H:c1}.
\begin{corollary}\label{H:c1}
Let $f(\x)=\sum\limits_{\ell=0}^\infty f_\ell(\x)$, where $f_\ell$
is a homogeneous polynomial of degree $\ell$, and
the convergence is uniform on $\s^k$.
Let $Y$ be a spherical $h$-harmonic
of degree $m$. Then
\begin{gather}\label{H:eq4}
\int_{\s^k} h^2(\x)f(\x)Y(\x)\,d\s(\x)
=\sum_{r=0}^\infty
\frac{1}{2^{m+2r-1}r!}\frac{\prod\limits_{i=0}^k \Gamma(\alpha_i+\frac12)}{\Gamma(m+\mu+r+1)}
\,\Delta_h^rY(\Dbold)f_{m+2r}(\x) .
\end{gather}
\end{corollary}

\section[Niven's formula for external ellipsoidal $h$-harmonics]{Niven's formula for external ellipsoidal $\boldsymbol{h}$-harmonics}\label{N}

We consider the ellipsoid $\J$ from \eqref{I:J} with semi-axes
$d_j:=\sqrt{t-a_j}$.
By substituting $y_j=d_jx_j$, $\y\in\J$ is transformed to $\x\in\s^k$.
If $f(\x)$ is a continuous function on $\s^k$, then
\begin{gather}\label{N:trans1}
 g(\y):=f\left(\frac{y_0}{d_0},\dots,\frac{y_k}{d_k}\right)
\end{gather}
becomes a continuous function on $\J$ and
\begin{gather}\label{N:trans2}
\int _{\s^k} f(\x)\,d\s(\x)=\prod_{j=0}^k d_j^{-1} \int_{\J} w(\y)g(\y)\,d\J(\y)
\end{gather}
with $w(\y)$ from \eqref{I:w}.

Let $Y$ denote a spherical $h$-harmonic of degree $m$, and
let $f$ be as in Corollary \ref{H:c1}.
We use the substitution $y_j=d_jx_j$ and \eqref{N:trans2} to write \eqref{H:eq4}
as
\begin{gather}
\int_\J h^2(\y)w(\y)g(\y)Y\left(\frac{y_0}{d_0},\dots,\frac{y_k}{d_k}\right)\,d \J(\y)
\nonumber\\
\qquad{}=\prod_{j=0}^k {d_j}^{2\alpha_j+1}\label{N:eq1}
\sum_{r=0}^\infty \frac{1}{2^{m+2r-1}r!}\frac{\prod\limits_{j=0}^k\Gamma(\alpha_j+\tfrac12)}
{\Gamma(m+\mu+r+1)}
\Lambda^rY(d_0\D_0,\dots,d_k\D_k)g_{m+2r}(\x),
\end{gather}
where $\Lambda$ is the operator
\[ \Lambda:=d_0^2\D_0^2+\dots+d_k^2\D_k^2 \]
and $g$ is def\/ined by \eqref{N:trans1}
with the same connection between $f_\ell$ and $g_\ell$.
We use \eqref{N:eq1} in the special case
$Y=G_\nn$, $m=2|\n|+|\p|$, and $g=\Phi(\cdot,\z)$, where $\|\z\|>d_0$
(the largest semi-axis of~$\J$).
We also note that
\[ E_\nn(t)G_\nn\left(\frac{y_0}{d_0},\dots,\frac{y_k}{d_k}\right)=F_\nn(\y),\]
see \cite[(6.1)]{Vo}.
Then the integral representation
\eqref{I:int2} for the external ellipsoidal $h$-harmo\-nic~$\F_\nn$
and \eqref{N:eq1} yield
\begin{gather*}
\F_\nn(\z)=\frac{2(\mu+m)}{E_\nn(t)}\sum_{r=0}^\infty \frac{1}{2^{m+2r-1}r!}
\frac{\prod\limits_{j=0}^k\Gamma(\alpha_j+\tfrac12)}
{\Gamma(m+\mu+r+1)}\\
\phantom{\F_\nn(\z)=}{} \times
\Lambda^rG_\nn(d_0\D_{x_0},\dots,d_k\D_{x_k})\Phi_{m+2r}(\x,\z).
\end{gather*}
We now use the identity
\[ G_\nn(d_0\D_0,\dots,d_k\D_k)g(\x)=E_\nn(t)G_\nn(\Dbold)g(\x) \]
which holds for any $h$-harmonic function $g$; see the formula before (6.3) in
\cite{Vo}. Then we replace the operator $-\Lambda$ by
$a_0\D_0^2+\dots+a_k\D_k^2$ which is possible because $\Phi(\cdot,\z)$ is $h$-harmonic.
Finally, we use Lemma \ref{F:l1} and arrive at the following theorem.
\begin{theorem}\label{N:t1}
The external ellipsoidal $h$-harmonic $\F_\nn$
admits the differential-difference  representation
\begin{gather}\label{N:Niven2}
\F_\nn(\z)=\sum_{r=0}^\infty \frac{(-1)^{r+m}(\mu+m)}{2^{m+2r}r!(\mu)_{m+r+1}}
\big(a_0\D_0^2+\dots+a_k\D_k^2\big)^r
G_\nn(\Dbold)\|\z\|^{-2\mu},
\end{gather}
where $m=2|\n|+|\p|$ and $\mu$ is defined in \eqref{1:mu}.
The expansion is valid for $\|\z\|^2>a_k-a_0$.
\end{theorem}
Formula \eqref{N:Niven2} is true for all $\z$ with $\|\z\|^2>a_k-a_0$ because we can choose
the number $t$ def\/ining the ellipsoid \eqref{I:J} very close to $a_k$.

We can rewrite \eqref{N:Niven2} slightly by using the connection between internal and
external spherical $h$-harmonics \cite[Corollary 2]{Vo}:
\[ \Y(\x)=\frac{(-1)^m}{2^m(\mu)_m} Y(\Dbold)\|\x\|^{-2\mu} .\]
This leads to the following corollary.
\begin{corollary}\label{N:c1}
The external ellipsoidal $h$-harmonic $\F_\nn$
admits the differential-difference representation
\begin{gather}\label{N:Niven3}
\F_\nn(\z)=\sum_{r=0}^\infty \frac{(-1)^r}{2^{2r}r!(\mu+m+1)_r}\big(a_0\D_0^2+\dots+a_k
\D_k^2\big)^r \G_\nn(\z)
\end{gather}
for $\|\z\|^2>a_k-a_0$.
\end{corollary}

\section{Examples}\label{P}
We illustrate some results of this paper in the planar case $k=1$.
We choose $a_0=-1$, $a_1=1$ and set $\alpha:=\alpha_1-\frac12$, $\beta:=\alpha_0-\frac12$.
If $p_0=p_1=0$ equation \eqref{E:Fuchs} becomes
\begin{gather}\label{Ex:Jacobi}
(t^2-1)\left(v''+\frac{\alpha+1}{t-1}v'+\frac{\beta+1}{t+1}v'\right)+\lambda_0 v=0.
\end{gather}
This equation is satisf\/ied by the Jacobi polynomial $v(t)=P_n^{(\alpha,\beta)}(t)$
when $\lambda_0=-n(n+\alpha+\beta+1)$.
Therefore, the Stieltjes polynomial $E_{n,\0}$ is a multiple of
$P_n^{(\alpha,\beta)}$. In \cite{Vo} the Stieltjes polynomials were normalized such that
their leading coef\/f\/icient equals $1$. This gives
\[ E_{n,\0}=a_n P_n^{(\alpha,\beta)},\]
where
\[ a_n=a_n^{(\alpha,\beta)}:=\frac{n!2^n}{(\alpha+\beta+n+1)_n} .\]
There are similar formulas for the other three parity vectors $\p$; compare \cite[\S~2]{Vo}.

A second solution of \eqref{Ex:Jacobi} with $\lambda_0=-n(n+\alpha+\beta+1)$
is $Q_n^{(\alpha,\beta)}$; see \cite[\S~4.61]{Sz}.
Comparing the behavior at inf\/inity, we see that
\[ \E_{n,\0}=b_n Q_n^{(\alpha,\beta)},\]
where
\[ b_n=b_n^{(\alpha,\beta)}:=2^{-n-\alpha-\beta}\frac{\Gamma(2n+\alpha+\beta+2)}
{\Gamma(n+\alpha+1)\Gamma(n+\beta+1)} .\]
Again, similar formulas hold for the other three parity vectors $\p$.

The planar sphero-conal coordinates $r,s_1$ as def\/ined in \cite[\S~4]{Vo}
essentially agree with polar coordinates $x_0=r\cos\phi$, $x_1=r\sin\phi$.
The connection is given by $\cos(2\phi)=s_1$.
The internal sphero-conal $h$-harmonic $G_{n,\0}$ is equal to
$r^{2n}E_{n,\0}(s_1)$ in sphero-conal coordinates.
Therefore, in cartesian coordinates, we obtain
\begin{gather*}
 G_{n,\0}(x_0,x_1)=a_n\big(x_0^2+x_1^2\big)^nP_n^{(\alpha,\beta)}
\left(\frac{x_0^2-x_1^2}{x_0^2+x_1^2}\right).
\end{gather*}
The external sphero-conal $h$-harmonic $\G_{n,\0}$ is given by
\begin{gather}\label{Ex:GG}
\G_{n,\0}(x_0,x_1)=a_n\big(x_0^2+x_1^2\big)^{-n-\alpha-\beta-1}P_n^{(\alpha,\beta)}
\left(\frac{x_0^2-x_1^2}{x_0^2+x_1^2}\right).
\end{gather}
The internal ellipsoidal $h$-harmonic $F_{n,\0}$ takes the form
$E_{n,\0}(t_0)E_{n,\0}(t_1)$ in ellipsoidal coordinates.
We can express $t_0$, $t_1$ explicitly as functions of $x_0$, $x_1$
as follows
\begin{gather}
t_0=\tfrac12\big(x_0^2+x_1^2\big)^2+\left(\tfrac14\big(x_0^2+x_1^2\big)^2-x_0^2+x_1^2+1\right)^{1/2},\label{Ex:t0}\\
t_1=\tfrac12\big(x_0^2+x_1^2\big)^2-\left(\tfrac14\big(x_0^2+x_1^2\big)^2-x_0^2+x_1^2+1\right)^{1/2}.\label{Ex:t1}
\end{gather}
It should be mentioned that analogous formulas do not exist when $k\ge 4$ and would be very
complicated when $k=2$ or $k=3$.
Now internal ellipsoidal $h$-harmonics can be written as
\begin{gather*}
F_{n,\0}(x_0,x_1)=a_n^2P_n^{(\alpha,\beta)}(t_0)P_n^{(\alpha,\beta)}(t_1)
\end{gather*}
with $t_0$, $t_1$ given in \eqref{Ex:t0}, \eqref{Ex:t1}, respectively.
Moreover, external ellipsoidal $h$-harmonics assume the form
\begin{gather}\label{Ex:FF}
\F_{n,\0}(x_0,x_1)=a_nb_n Q_n^{(\alpha,\beta)}(t_0)P_n^{(\alpha,\beta)}(t_1).
\end{gather}
Formula \eqref{Ex:FF} is meaningful for all $(x_0,x_1)$ which do not lie on the common focal line
$-\sqrt2\le x_0\le\sqrt2$, $x_1=0$ of the confocal ellipses determined by a constant value
of the ellipsoidal coordinate $t_0$.

Formula \eqref{N:Niven3} now expresses the function \eqref{Ex:FF} in terms of the
function \eqref{Ex:GG} as follows.
\begin{theorem}\label{Ex:thm1}
Let $n\in\N_0$. For every $(x_0,x_1)\in\R^2$ with $x_0^2+x_1^2>2$ we have
\begin{gather}
b_nQ_n^{(\alpha,\beta)}(t_0)P_n^{(\alpha,\beta)}(t_1)
\nonumber\\
\qquad{}=\sum_{r=0}^\infty \frac{(-1)^r}{2^{2r}r!(2n+\alpha+\beta+2)_r}
\big(\D_1^2-\D_0^2\big)^r
\big(x_0^2+x_1^2\big)^{-n-\alpha-\beta-1}P_n^{(\alpha,\beta)}
\left(\frac{x_0^2-x_1^2}{x_0^2+x_1^2}\right),\!\!\label{Ex:Niven}
\end{gather}
where on the left-hand side $t_0$, $t_1$ are given by
\eqref{Ex:t0}, \eqref{Ex:t1}, respectively.
\end{theorem}
Obviously, we could replace $\D_1^2-\D_0^2$ by $2\D_1^2$ on the right-hand side of \eqref{Ex:Niven}.
Formula \eqref{Ex:Niven} expands the left-hand side in a series of homogeneous functions
whose degree decreases with $r$. The convergence will be especially good when $x_0^2+x_1^2$
is large. The formula can be checked with symbolic computer algebra by
replacing $x_j$ by $ux_j$ on the left-hand side and then expanding in an asymptotic
series at $u=\infty$. Computer algebra conf\/irmed the correctness of \eqref{Ex:Niven}
for small values of $n$. The natural question arises whether
this formula or the more general formula \eqref{N:Niven3} can be proved more
directly, that is, without using any integrals in the proof.

\pdfbookmark[1]{References}{ref}
\LastPageEnding
\end{document}